\documentclass[a4paper,11pt]{article}
\usepackage{a4wide}
\usepackage{amssymb}
\usepackage{amsmath}
\usepackage{enumerate}
\usepackage{amsthm}
\usepackage{amsfonts}

\newtheorem{thm}{Theorem}[section]
\newtheorem{lem}{Lemma}[section]
\newtheorem{df}{Definition}[section]
\newtheorem{prop}{Proposition}[section]

\newtheorem{cor}{Corollary}[section]

\newcommand{\C}{\mathbb{C}}
\newcommand{\R}{\mathbb{R}}
\newcommand{\N}{\mathbb{N}}

\newcommand{\DD}{\mathcal{D}}
\newcommand{\FF}{\mathcal{F}}

\newcommand{\e}{\varepsilon}

\def\limess{\qopname\relax m{lim\,ess}}
\DeclareMathOperator*{\esssup}{ess\,sup}
\DeclareMathOperator*{\essinf}{ess\,inf}
\DeclareMathOperator*{\sgn}{sgn}

\begin{document}

\title{Large time behaviour of a conservation law regularised by a Riesz-Feller operator: the sub-critical case}

%    Remove any unused author tags.

%    author one information
%\author{C. M. Cuesta}
%\address{Dpto. Matem\'aticas, Facultad de Ciencia y Tecnolog\'{i}a (UPV/EHU), Barrio Sarriena s/n 48940 Leioa-Bizkaia, Spain}
%\curraddr{}
%\email{carlotamaria.cuesta@ehu.eus}
%\thanks{}

\author{C.~M. Cuesta, X. Diez-Izagirre\\
  University of the Basque Country (UPV/EHU),\\  Faculty of Science and Technology, Department of Mathematics, \\Barrio Sarriena S/N, 48940 Leioa, Spain.\\
carlotamaria.cuesta@ehu.eus, xuban.diez@ehu.eus}

\date{June 6, 2022}%Received: date / Accepted: date}

\maketitle

\begin{abstract}
We study the large time behaviour of the solutions of a non-local regularisation of a scalar conservation law. This regularisation is given by a fractional derivative of order $1+\alpha$, with $\alpha\in(0,1)$, which is a Riesz-Feller operator. The non-linear flux is given by the locally Lipschitz function $|u|^{q-1}u/q$ for $q>1$. We show that in the sub-critical case, $1<q < 1 +\alpha$, the large time behaviour is governed by the unique entropy solution of the scalar conservation law. Our proof adapts the proofs of the analogous results for the local case (where the regularisation is the Laplacian) and, more closely, the ones for the regularisation given by the fractional Laplacian with order larger than one, see Ignat and Stan (2018). The main difference is that our operator is not symmetric and its Fourier symbol is not real. We can also adapt the proof and obtain similar results for general Riesz-Feller operators. 
\end{abstract}


%-----------------------------------------------------------------------------------------------------------------------------------------------------------------------------------------------------------------

\section{Introduction}\label{introduction}
In this paper, we study the large time asymptotic behaviour of non-negative solutions to the convection-diffusion equation
\begin{equation}\label{NCL}
\begin{cases}
\partial_t u(t,x) + |u(t,x)|^{q-1}\partial_x u(t,x) = \partial_x \DD^\alpha[u(t,\cdot)](x), \quad &t>0, \ \ x\in \R, \\
u(0,x) = u_0(x), \quad &x\in \R,
\end{cases}
\end{equation}
where $u_0 \in L^1(\R)\cap L^\infty(\R)$ and $q>1$. The operator, acting on $x$, $\DD^\alpha[\cdot]$ has $\alpha\in(0,1)$ and is defined by means of
\begin{equation}\label{WMFD}
\DD^\alpha[g](x) = d_{\alpha+1} \int_{-\infty}^{0}\frac{g(x+z)-g(x)}{|z|^{\alpha+1}} dz, \quad \text{for} \ \ 0<\alpha<1, \ \ d_{\alpha+1} = \frac{1}{\Gamma(-\alpha)}.
\end{equation}
The operator $\partial_x \DD^\alpha[\cdot]$ is of Riesz-Feller type, as we shall see below. The operator $\DD^\alpha[\cdot]$ can also be seen as a right Weyl-Marchaud fractional derivative (see \cite{Marchaud}, \cite{Weyl}) of order $\alpha$. The non-linear flux $f(u) = |u|^{q-1} u/q$, is considered here as a paradigm locally Lipschitz function.

The equation in (\ref{NCL}) is a modified Burgers' equation, and appears in \cite{SK1985} as a model of viscoelastic waves with $\alpha=1/2$. There are other models of physical phenomena where this kind of non-local operator appears, such as problems in fluid dynamics, see for instance the references listed in \cite{AU}.

The models that motivate our study describe the internal structure of hydraulic jumps in a shallow water model. The general form being, for $C_1,C_2 \geq 0$,
\[
\partial_t u + \partial_x(f(u)) =C_1 \partial_x \DD^{1/3}[u] +C_2 \partial_x^3 u, \quad t>0, \quad x\in \R
\]
see \cite{Kluwick2010} and \cite{NV}, where the flux might be $u^2$ or $u^3$ and the dispersive term might or might not be relevant, depending on the asymptotic regime considered.

Another example where such operators appear (although not as a regularising term) can be found in \cite{Fowler}, where a model for dune formation is presented.

Formally, the study of the large time behaviour can be transferred to a limit problem by the appropriate scaling; for any $\lambda>0$, let the change of variables
\begin{equation}\label{lambda:scaling}
 t= \lambda^q s \quad x = \lambda y
\end{equation}
and the function 
\begin{equation}\label{u:scaling}
u_\lambda(s,y) := \lambda \, u(\lambda^q s, \lambda y).
\end{equation}
Then, if $u$ is a solution of (\ref{NCL}), $u_\lambda$ satisfies
\begin{equation}\label{lambda:Prob}
\begin{cases}
\partial_s u_\lambda + |u_\lambda|^{q-1} \partial_y u_\lambda = \lambda^{q-1-\alpha} \partial_y \DD^\alpha\left[u_\lambda(s,\cdot)\right](y),& \quad s>0, \ y\in \R, \\
u_\lambda(0,y) = \lambda u_0(\lambda y),& \quad y\in \R.
\end{cases}
\end{equation}
Observe that when $t\to\infty$, if we keep $s$ of order one, then $\lambda\to \infty$, and in the new variables this means that, depending on the sign of the exponent $q-1-\alpha$, different terms dominate the limit behaviour. According to this heuristic argument, we can distinguish three different regimes: The sub-critical case $1<q<1+\alpha$ (the conservation law formally dominates), the critical case $q= 1+\alpha$ (one expects self-similar behaviour associated to the balance of all terms) and the super-critical case $1+\alpha<q$ (the non-local heat equation formally dominates).

In this paper we focus on the sub-critical case for non-negative solutions, henceforth we assume that $1<q<1+\alpha$. Our main theorem is:
\begin{thm}\label{subcritical:asymptotics}
  Let $1<q<1+ \alpha$ and $1\leq p < \infty$. Given $u_0 \in L^\infty(\R)\cap L^1(\R)$ with $\int_{\R}u_0(x)dx=M>0$ and $u_0(x)\geq 0$ for all $x\in\R$, then, there exists a unique solution (\ref{NCL}), with,
  \begin{itemize}\item[(i)] for $p=1$, $\partial_t u \in C((0,\infty), L^1(\R))$, $u \in C((0,\infty), C_b(\R) \cap C^1(\R))$ and $u \in C([0,\infty), L^1(\R))$,
\item[(ii)] for $1< p < \infty$, $\partial_t u \in C((0,\infty), L^p(\R))$ and $u \in C((0,\infty), L^p(\R) \cap\dot{H}^{\theta,p}(\R))$ for any $0<\theta< 1+\alpha + \min\{\alpha, q-1\}$.
  \end{itemize}
Moreover, $u$ satisfies 
\begin{equation}\label{asymp:sub}
\lim_{t\to \infty} t^{\frac{1}{q}(1 - \frac{1}{p})} \|u(t,\cdot) - U_M(t,\cdot)\|_{L^p(\R)}=0,
\end{equation}
where $U_M$ is the unique entropy solution of %the purely convective equation
\begin{equation}\label{inviscid:Prob}
\begin{cases}
\partial_t U_M + \partial_x ( |U_M|^{q-1} \, U_M /q ) = 0,& \quad  t>0, \ x\in \R,\\
U_M(0,x) = M \delta_0,& \quad x\in \R.
\end{cases}
\end{equation}
\end{thm}
Here the notation $\dot{H}^{\theta,p}(\R)$ is for the homogeneous Sobolev spaces (e.g. \cite{Adams}) of order $\theta>0$
\begin{equation}
\label{H:sp}
%  %L^p(R) \cap
  \dot{H}^{\theta,p}(\R):=\left\{ g \in \mathcal{S}'(\R) : \FF^{-1} \left[|\xi|^\theta \FF(g)\right] \in L^p(\R)\right\}
\end{equation}
and $\FF$ denotes the Fourier transform
\begin{equation}\label{Fourier:def}
\FF(g(x))(\xi)= \frac{1}{\sqrt{2\pi}} \int_{\R}{g(x) e^{-i\xi x} \ dx}.
\end{equation}

We also recall here the definition of entropy solution of (\ref{inviscid:Prob})
(see \cite{Kruzhkov}):
%..... SAY SMTH ABOUT UM... HERE?%In this section, we recall some classical result for the purely convective problem (\ref{inviscid:Prob}) and for (\ref{NCL}). For the former we recall the results of \cite{LP}, first we give the definition of entropy solution in the sense of Kru\v{z}kov's (see \cite{Kruzhkov}):
\begin{df}\label{entropy:SCL}
Let $U_M$ be a weak solution of (\ref{inviscid:Prob}) such that
\[
U_M \in L^\infty((0,\infty), L^1(\R)) \cap L^\infty((\tau,\infty)\times \R), \ \forall \tau \in (0,\infty),
\]
then $U_M$ is said to be an entropy solution of (\ref{inviscid:Prob}) if, and only if, the following inequality holds for every $k\in \R$ and $\varphi \in C_c^\infty((0,\infty) \times \R)$ non-negative
\begin{equation}\label{inviscid:entropy:inequality}
\int_{0}^{\infty}\int_{\R} \left( |U_M - k| \partial_t \varphi + \sgn(U_M-k) (f(U_M) - f(k)) \partial_x \varphi \right) \, dx dt \geq 0,
\end{equation}
with $f(u) = |u|^{q-1} u/q$, and for any $\psi \in C_b(\R)$
\begin{equation}\label{cons:mass:U:2}
\limess_{t \downarrow 0} \int_{\R} U_M(t,x) \psi(x) \, dx = M \psi(0).
\end{equation}
\end{df}
We recall that this unique entropy solution is given by the $N$-wave profile (see \cite{LP}).

The proof follows the method developed by Kamin and V\'azquez in \cite{KV}, namely, that with the rescaling (\ref{lambda:scaling})-(\ref{u:scaling}), (\ref{asymp:sub}) is formally equivalent to
\begin{equation}\label{asymp:ulambda}
  \lim_{\lambda\to\infty} \|u_\lambda(s_0,\cdot) - U_M(s_0,\cdot)\|_{L^p(\R)}=0
  %\rightarrow 0, \quad \text{as} \quad \lambda \to \infty,
\end{equation} 
for some $s_0>0$ fixed. We recall that, the local case in $\R^N$ for $N\in \N$ and for all $q>1$ has been analysed by Escobedo, V\'azquez and Zuazua \cite{EVZ1, EVZ2, EZ}. %In their model, the non-linear propagation term takes place in a given direction. The results for the critical and supercritical case are proved in \cite{EZ}, and the ones for the subcritical case in \cite{EVZ1, EVZ2}.

The sub-critical case, for the fractional Laplacian as viscous term, in dimension one has been studied by Ignat and Stan in \cite{IS}, and these are the results that we shall adapt. We also mention that, Biler, Karch and Woyczy\'{n}ski in \cite{BKW1,BKW2} study the critical and super-critical cases for a more general L\'evy  operator as viscous term, showing the expected asymptotic behaviour; this being given by the self-similar solution and by the fractional heat kernel, respectively. This operator is defined by means of a Fourier multiplier, such that the symbol can be represented by the L\'evy-Khintchine formula in the Fourier variable (see \cite[Chapter~1, Theorem~1]{Bertoin}). The fractional Laplacian as a particular example of such operators. In all these results the non-negativity of the symbol and the symmetry of the operator are necessary conditions. The operator (\ref{WMFD}) does not satisfy these conditions. 

For completeness, we also consider the case of a general Riesz-Feller operator. Thus we show how the analogous to Theorem~\ref{subcritical:asymptotics} follows if the term $\partial_x \DD^\alpha[u]$ is replaced by the Riesz-Feller operator term $D^\beta_\gamma[u]$ of order $\beta\in (1,2)$ and skewness $\gamma$. We recall that such operator can be defined by means of a Fourier multiplier operator (see e.g. \cite{MLP}), we choose the definition
\begin{equation}\label{Riesz-Feller:operator}
\FF(D^\beta_\gamma[g])(\xi) = \psi^\beta_\gamma(\xi) \, \FF(g)(\xi),
\end{equation}
where $\beta \in (0,2]$ and $|\gamma| \leq \min\{\beta, 2-\beta\}$ and the symbol satisfies
\begin{equation}\label{Riesz-Feller:symbol}
\psi^\beta_\gamma(\xi) = - |\xi|^\beta \, e^{\displaystyle- i \sgn(\xi) \gamma \frac{\pi}{2}}.
\end{equation}
In particular, we observe that the derivative of (\ref{WMFD}) is of this form with $\beta= 1+\alpha$ and $\gamma = 1-\alpha$, as we explain in the next section.
We recall that the definition we use here for Riesz-Feller operators differs from the usual one: the symbol we obtain is the complex conjugate of the one with the standard definition. This is simply because such definition uses the complex conjugate of (\ref{Fourier:def}) (up to a scaling factor) as Fourier transform.

%NOTE: this is with our def of Fourier transform, it is just complex conjugation

Finally, we mention that we focus on non-negative solutions, that might be considered as certain density functions. Results for sign changing solutions for the local case and the one regularised by the fractional Laplacian appear in \cite[\S~3]{EVZ1} and \cite[\S~6]{CIP}, respectively.

The paper is organised as follows. Section~\ref{sec:2} contains preliminary results and is divided into three sections.

In the first one, we recall some properties of the non-local operator (\ref{WMFD}) and its derivative. These include some equivalent integral representations, the computation of their Fourier symbol and the definition and properties of certain {\it dual} operators. In the second part we recall the linear problem and derive some estimates of the fundamental solutions that will be necessary later. We end the section define mild solutions for problem (\ref{NCL}) and give existence and regularity results for this. We also recall the viscous entropy inequality, as well as derive a comparison principle. Most of the results appear in \cite{DC} or can be adapted from \cite{IS} with the ingredients given in the previous sections, thus we shall only prove which does not appear elsewhere.

In Section~\ref{sec:3} we derive necessary {\it a priori} estimates, namely a Oleinik type entropy inequality and an energy type estimate. In order to do this we first consider the problem with a positive initial condition (which makes the non-linear flux regular, since positivity is preserved, by the comparison principle). In this case we can show the desired estimates, which are then preserved in the limit to a non-negative initial condition. The proofs are similar to those in \cite{IS}, we only give the details where necessary.

In Section~\ref{sec:4} we prove Theorem~\ref{subcritical:asymptotics}. First, we translate the results of Section~\ref{sec:3} into the rescaled problem (\ref{lambda:Prob}). Additionally, before proving the limit $\lambda\to\infty$, we have to take care of the behaviour of $u_\lambda$ for large $|y|$. With these estimates we prove the limit (\ref{asymp:ulambda}).

We finish in Section~\ref{sec:5} by showing how the proofs generalise to Riesz-Feller operators.

We recall that we do not give the proofs that are analogous to those given in \cite{IS}, the interested reader is also referred to Diez-Izagirre's PhD thesis \cite[Chapter~3]{XDThesis} for details.

%---------------------------------------------------------------------------------------------------------------------------------------------------------------------

\section{Preliminary results}\label{sec:2}%{preliminary:results}
In this section we recast the necessary Lemma's that are needed to adapt the results of \cite{IS}.

\subsection{Derivation and integration by parts rules}% of non-local operators}
%In this section we give alternative formulations, which under certain smoothness conditions are equivalent, of the operator (\ref{WMFD}) and its derivative.
In this section we derive integration by parts rules for (\ref{WMFD}).
%We also recall an integration by parts rule and product rule for the fractional Laplacian (see  \cite{}).

First we need the notation for the following operator:
\begin{equation}
\label{leftWM:alpha}
\overline{\DD^\alpha}[g](x)= - d_{\alpha+1} \int_{0}^{\infty} \frac{g(x+z)-g(x)}{|z|^{\alpha+1}}\ dz.
\end{equation}

%The following result can be found in \cite{DC}.
%\begin{lem}[Equivalent representations of $\partial_x \DD^\alpha$ and $\DD^\alpha$]\label{equiv:repre}
% If $\alpha \in (0,1)$, then for all $\varphi \in C_b^1(\R)$ and all $x\in \R$,
%\[
%\DD^\alpha [\varphi] (x)= d_\alpha \int_{-\infty}^{x} \frac{\varphi'(z)}{(x-z)^\alpha} \, dz.
%%\label{integralformula2}
%\]

%Moreover, for all $\varphi \in C^2_b(\R)$ and all $x\in \R$,
%\[
%\begin{split}
%\partial_x \DD^\alpha[\varphi](x) &= d_{\alpha+2} \int_{-\infty}^{0} \frac{\varphi(x+z)-\varphi(x)- \varphi'(x)z}{|z|^{\alpha+2}} dz \\
%&= d_\alpha \int_{-\infty}^{x} \frac{\varphi''(z)}{(x-z)^\alpha}\, dz.
%%= {}^C\DD^{\alpha+1}[\varphi](x).
%\end{split}
%%\label{integralformula1}
%\]
%\end{lem}

Then the following integration by parts follows (see \cite{DC}).
\begin{lem}\label{integrationbyparts}%INTEGRATING IN TIME NOT NECESSARY HERE 
Let $\alpha\in(0,1)$, $u\in C^2_b(\R)$ and $\varphi\in C^\infty_c(\R)$. Then, for all $t>0$
\begin{equation*}
  \int_{\R} \varphi(x) \partial_x \DD^\alpha[u(t,\cdot)](x) \, dx =
  \int_{\R} \partial_x \overline{\DD^\alpha}[\varphi](x) u(t,x) \, dx,
\end{equation*}
\end{lem}

We also recall the following representation of the operators (see also \cite{DC}).
\begin{lem}[Integral representations]\label{equiv:repre}
If $\alpha \in (0,1)$, then for all $g \in C_b^1(\R)$ and all $x\in \R$,
\[
\partial_x \DD^\alpha[g](x) = d_{\alpha+2} \int_{-\infty}^{0} \frac{g(x+z)-g(x)- g'(x)z}{|z|^{\alpha+2}} dz,
\]
and
\[
\partial_x \overline{\DD^\alpha}[g](x)  = d_{\alpha+2} \int_{0}^{\infty} \frac{g(x+z)-g(x)- g'(x)z}{|z|^{\alpha+2}}\ dz.
\]
\end{lem}

With (\ref{Fourier:def})
%definition and notation for the Fourier transform,%DO WE USE THE HAT NOTATION??? NOT REALLY
%\begin{equation}\label{Fourier:def}
%\FF(g(x))(\xi)= \frac{1}{\sqrt{2\pi}} \int_{\R}{g(x) e^{-i\xi x} \ dx},
%\end{equation}
we obtain, formally, the Fourier symbol of $\DD^\alpha[\cdot]$ (see e.g. \cite[Chapter~7]{SKM}),
\begin{equation}
\FF\left( \DD^\alpha[g](x) \right)(\xi) = (i \xi)^\alpha  \FF(g)(\xi),%, \quad \text{for} \ 0<\alpha<1,
\label{f:symbol:D:alpha}
\end{equation}
and that of $\overline{\DD^\alpha}[\cdot]$,
\begin{equation}
\FF\left( \overline{\DD^\alpha}[g](x) \right)(\xi) = - (-i \xi)^\alpha  \FF(g)(\xi).%, \quad \text{for} \ 0<\alpha<1.
\label{f:symbol:adjoint:D:alpha}
\end{equation}

In particular $\DD^\alpha[\cdot]$ is not of Riesz-Feller type because its symbol has $\beta=\alpha$ and $\gamma=2-\alpha$, but $\overline{\DD^\alpha}[\cdot]$ belongs to this class, since $\beta=\alpha$ and $\gamma=\alpha$. We also observe that their symbols satisfy $(i\xi)^\alpha=- \overline{(-(-i\xi)^\alpha)}$, 
where the bar on the right-hand side denotes complex conjugation.
%KEEP THIS FOR CONSISTENCY; GIVE EQ. NUMBER?; THE REASON OF THE DEFINITION IT THE INTEGRATION BY PARTS FORMULA, ()DISPLAY

We can then conclude that
\[
%\label{dx:D:symbol}
\FF(\partial_x \DD^{\alpha}[g])(\xi)=(i\xi)^{\alpha+1} \, \FF(g)(\xi).
\]
We observe that, writing
\[
(i\xi)^{\alpha+1} = -|\xi|^{\alpha+1} \left( \cos\left( (1-\alpha) \frac{\pi}{2} \right) -i \sgn(\xi) \sin\left( (1-\alpha) \frac{\pi}{2} \right) \right),
\]
$\partial_x \DD^{\alpha}[\cdot]$ is an operator of Riesz-Feller type with $\beta=1+\alpha$ and $\gamma=1-\alpha$ (see (\ref{Riesz-Feller:operator})-(\ref{Riesz-Feller:symbol})). %We recall that the definition we use here for Riesz-Feller operators differs from the usual one: the symbol we obtain is the complex conjugate of the one with the standard definition. This is simply because such definition uses the complex conjugate of (\ref{Fourier:def}) (up to a scaling factor) as Fourier transform.

We also then get:%HERE HERE: DO WE NEED THIS? DONT KNOW
\[
%\label{dx:Dadjoint:symbol}
\FF(\partial_x\overline{\DD^{\alpha}}[g])(\xi)=(-i\xi)^{\alpha+1} \ \FF(g)(\xi), \quad \text{for} \ 0<\alpha<1.
\]

%DO WE NEED THIS?? NO
%Apart from this, while the Caputo type operator is bounded from $C_b^1$ to $C_b$ and from $H^{m+\alpha}$ to $H^m$, it is clear that (\ref{WMFD}) exists for a far more general assumption such as bounded functions satisfying local H\"older condition for some $\gamma>\alpha$. See \cite[Chapter~2]{SKM} for more insight into fractional derivatives.
%........

With this Fourier representation formulas we can now prove the following integration by parts rule:% for sufficiently regular functions:
\begin{lem}\label{partial:int:parts:D:alpha}
  Let $0<\alpha<1$ and $1/2<\theta_1$, $\theta_2<1$ such that $1+\alpha = \theta_1 + \theta_2$, assume also that $g\in H^2(\R)$ so that $\partial_x\DD^\alpha[g]$, $\DD^{\theta_1}[g]$, $\overline{\DD^{\theta_2}}[h] \in L^2(\R)$, and let $h \in  L^2(\R)$, then 
\[
\int_{\R} \partial_x \DD^\alpha[g](x) \, h(x) \, dx = - \int_{\R} \DD^{\theta_1}[g](x) \, \overline{\DD^{\theta_2}}[h](x) \, dx.
\]
\end{lem}
\proof
Since $\partial_x \DD^\alpha[g]$, $h \in L^2(\R)$, then Plancherel's theorem yields
\[
\begin{split}
  \int_{\R} \partial_x \DD^\alpha[g](x) \, h(x) \, dx &= \int_{\R} (i \xi)^{1+ \alpha} \, \FF(g)(\xi) \overline{\FF(h)(\xi)} \, d\xi \\
&  = -\int_{\R} (i \xi)^{\theta_1} \, \FF(g)(\xi) \, \overline{(-(-i \xi)^{\theta_2})\FF(h)(\xi)} \, d\xi \\
&  =- \int_{\R}  \FF(\DD^{\theta_1}[g])(\xi) \, \overline{\FF(\overline{\DD^{\theta_2}}[h])(\xi)} \, d\xi \\ & = 
% - \int_{\R} \DD^{s_1}[g](x) \, \overline{\overline{\DD^{s_2}}[h](x)} \, dx =
- \int_{\R} \DD^{\theta_1}[g](x) \, \overline{\DD^{\theta_2}}[h](x) \, dx,
\end{split}
\]
for $1/2<\theta_1$, $\theta_2<1$ such that $1+\alpha = \theta_1 + \theta_2$.
%, where we have used
%%the fact that $g(x)$ and $\overline{\DD^{s_2}}[h](x)$ are real functions and
%the identity $(i \xi)^{s_2} = \overline{(- i \xi)^{s_2}}$ (the bar indicating complex conjugation).%% in order to rewrite the integral term. %Notice that while the line over the operator denotes the adjoint of the operator, the line over denotes complex conjugation. 
\endproof

We now recall some facts about the fractional Laplacian that we need later
to conclude $L^p$-regularity of the solution.% and an energy estimate???. 

There are several equivalent definitions of the fractional Laplacian (see \cite{Kwasnicki}). Here we consider the one given by the Fourier symbol for $0<\theta<2$, namely,
\begin{equation}\label{Fractional:Laplacian}
|D|^\theta[g](x): = (-\Delta)^{\theta/2}[g](x):= \FF^{-1}\left(|\xi|^\theta \FF(g)(\xi)\right)(x).
\end{equation}
We observe that applying Plancherel's theorem $\||D|^\theta[g]\|_{L^2(\R)}=\|g\|_{\dot{H}^{\theta,2}(\R)}$ and also
\begin{equation}
  \label{homo:norm}
  \|\DD^\theta[g]\|_{L^2(\R)}=\|\overline{\DD^\theta}[g]\|_{L^2(\R)}=\|g\|_{\dot{H}^{\theta,2}(\R)} \quad \mbox{for}\quad 0<\theta<1
\end{equation}
since
\[
\|\DD^\theta[g]\|_{L^2(\R)} = \| (i \, \cdot)^\theta \FF(g)(\cdot) \|_{L^2(\R)} = \| |\cdot|^\theta \FF(g)(\cdot)\|_{L^2(\R)} = \||D|^\theta[g]\|_{L^2(\R)}.
\]
For comparison purposes, we give an equivalent definition of the fractional Laplacian of order $\theta>1$ for functions $g\in C_b^2(\R)$ (see \cite{Droniou3}):
\begin{equation}\label{int:frac:Lap}
(-\Delta)^{\theta/2}[g]=c_\theta\int_{\R} \frac{g(x+z)-g(x)-g'(x)z}{|z|^{\theta+1}}dz\quad \mbox{with}\quad \theta \in (1,2)
\end{equation} with $c_\theta=\frac{2^\theta\Gamma((\theta+1)/2)}{\pi^{\frac{1}{2}}\Gamma(-\theta/2)}$.

\subsection{Linear fractional diffusion equation}
In this section we recall some results concerning the linear problem
\begin{equation}\label{Linear:Prob}
\begin{cases}\displaystyle
\partial_t U(t,x) - \partial_x \DD^\alpha\left[U(t,\cdot)\right](x)=0, \quad &t>0, \ x\in \R, \\
U(0,x)= u_0(x), \quad &x\in \R,
\end{cases}
\end{equation}
for initial data $u_0 \in L^\infty(\R)$ the solution of (\ref{Linear:Prob}) can be represented as follows,
\[
U(t,x)= \left(K(t,\cdot)*u_0\right)(x) = \int_{\R}{K(t,x-y)\, u_0(y) \, dy},
\]
such that the kernel, $K(t,x)$, is defined by means of
\begin{equation}
\label{kernel:CP}
K(t,x) = \FF^{-1}\left( e^{(i \xi)^{\alpha+1} t}\right)(x), \quad \forall t>0, \ x\in \R
\end{equation}
which can be formally obtained using Fourier transform (see \cite{AK2} for the proof). Some pertinent properties of the kernel are derived in \cite{AHS} and \cite{DC}. In particular we recall that,
\begin{equation}\label{K:self:sim}
K(t,x) = \frac{1}{t^{\frac{1}{1+\alpha}}} K\left(1, \frac{x}{t^{\frac{1}{1+\alpha}}}\right), \quad \forall t>0,\ x\in \R,
\end{equation}
$K$ is non-negative and $K(t,x) \in C^\infty((0,\infty)\times \R)$, it also preserves mass and has the semi-group property.

Since the regularity of solutions is established with respect to derivation with the fractional Laplacian, we need the following estimates:
%(that will combine with Lemma~\ref{Chain:Product:rule} above for the non-linear problem).

\begin{lem}[Time behaviour of $K$]\label{K:time:behaviour} 
For all $\alpha, \theta \in (0,1)$, and $1\leq p \leq \infty$, $K(t,x)$ satisfies the following estimates for any $t>0$:
\[
\begin{split}
\|K(t,\cdot)\|_{L^p(\R)}& = C t^{-\frac{1}{1+\alpha} \left(1-\frac{1}{p} \right)}, \\ 
\| \partial_x K(t,\cdot)\|_{L^p(\R)}& \lesssim  t^{-\frac{1}{1+\alpha} \left(1-\frac{1}{p}\right)- \frac{1}{1+\alpha}} ,\\
\| |D|^\theta [K(t,\cdot)]\|_{L^p(\R)}& \lesssim  t^{-\frac{1}{1+\alpha} \left(1-\frac{1}{p}\right)- \frac{\theta}{1+\alpha}}, \\ 
\| |D|^\theta [\partial_x K(t,\cdot)]\|_{L^p(\R)}& \lesssim  t^{-\frac{1}{1+\alpha} \left(1-\frac{1}{p}\right)- \frac{1+\theta}{1+\alpha}}, \\
\end{split}
\] 
for some constant $C>0$.
\end{lem}
\proof
The first and second identities follow from (\ref{K:self:sim}), the mass-conservation property of $K$ and $\partial_x K$ and that they are bounded on $(0,T)\times\R$ for any $T>0$. 

For the third estimate we first use (\ref{K:self:sim}), then we rescale the fractional Laplacian:
\begin{equation}\label{D^sK:estimate:t}
\begin{split}
| |D|^\theta[K(t,\cdot)](x)| =& \frac{1}{t^{\frac{1}{1+\alpha}}} \left| |D|^\theta\left[K\left(1,\frac{\cdot}{t^{\frac{1}{1+\alpha}}}\right)\right](x)\right| \\  =&  \frac{1}{t^{\frac{1+\theta}{1+\alpha}}} \left| |D|^\theta[K(1,\cdot )]\left(\frac{x}{t^{\frac{1}{1+\alpha}}}\right)\right|.
\end{split}
\end{equation}
Now, when computing the $L^p$-norm we apply the change of variable $X= \frac{x}{t^{\frac{1}{1+\alpha}}}$:
\begin{equation}\label{estimate:D^sK}
\begin{split}
\|  |D|^\theta[K(t,\cdot)](x) \|_{L^p(\R)} =& \frac{1}{t^{\frac{1+\theta}{1+\alpha}}} \left(\int_{\R} \left| |D|^\theta[K(1,\cdot )]\left(\frac{x}{t^{\frac{1}{1+\alpha}}}\right)\right|^p \, dx \right)^{1/p} \\
=& \frac{1}{t^{\frac{1+\theta}{1+\alpha}}} t^{\frac{1/p}{1+\alpha}}  \left(\int_{\R} \left| |D|^\theta[K(1,\cdot )]\left( X \right)\right|^p \, dX \right)^{1/p}.
% \\ 
%\lesssim&  t^{-\frac{1}{1+\alpha} \left(1-\frac{1}{p}\right)- \frac{s}{1+\alpha}}.
\end{split}
\end{equation}
It only remains to prove that the $L^p$-norm of $|D|^\theta[ K(1,\cdot)]$ is finite. In order to show this, we first observe that using (\ref{Fractional:Laplacian}), the integrand of (\ref{estimate:D^sK}) is bounded:
\begin{equation}\label{D^s:K}
\begin{split}
\left| |D|^\theta[K(1,\cdot )]\left( X \right)\right| =& \frac{1}{\sqrt{2\pi}} \left| \int_{\R} |\xi|^\theta \, e^{(i \xi)^{1+\alpha}} e^{i X \xi} d\xi \right| \\
 \leq& \frac{1}{\sqrt{2\pi}} \int_{\R} |\xi|^\theta \, e^{-|\xi|^{1+\alpha} \sin\left( \frac{\alpha \pi}{2}\right)} \, d\xi < \infty.
\end{split}
\end{equation}
Next we show that $||D|^\theta[K(1,\cdot )(X)|^p$ is integrable for large $|X|$. We will apply \cite[Lemma~2]{PT}, which gives that for large $X$
  \[
  \int_0^\infty e^{-i X \xi}e^{-\sigma \xi^r}\xi^\theta d\xi=O(|X|^{-1-\theta})
  \]
  with $r\in(0,2)$, $\theta\geq 0$ and $\sigma$ such that
  \begin{equation}\label{cond:sigma}
\sigma\in \left\{ a+ i b \in \C : \ -\cos\left( \frac{r\pi}{2} \right) \leq a \leq 1, \ |b| \leq -\tan\left( \frac{r\pi}{2} \right) \right\}.
\end{equation}

  %of the complex plane. Namely that  We shall adopt the notation 
%\[
%\sigma= \sin\left(\frac{\alpha \pi}{2}\right) - i \cos\left(\frac{\alpha \pi}{2}\right).
%\]

%in order to apply \cite[Lemma~2]{PT}, which gives the estimates of for large $X$ of integrals of the form $\int_0^\infty e^{-iX\xi}e^{-\sigma \xi^r}\xi^sd\xi$ with $r\in(0,2)$, $s\geq 0$ and $\sigma$ in certain sector of the complex plane.

We first write,
  \[
\begin{split}
|D|^\theta[K(1,\cdot)](X) =& \frac{1}{\sqrt{2 \pi}}  \int_{\R} |\xi|^\theta \, e^{(i\xi)^{1+\alpha}} e^{i X \xi} \, d\xi \\
=&  \frac{1}{\sqrt{2 \pi}}  \int_{\R} |\xi|^\theta \, e^{-|\xi|^{1+\alpha} \left(\sin\left(\frac{\alpha \pi}{2}\right) - i \sgn(\xi)\cos\left(\frac{\alpha \pi}{2}\right)  \right)} e^{i X \xi} \, d\xi \\
%=& \frac{1}{\sqrt{2 \pi}}  \int_{0}^{\infty}  \xi^s \, e^{-\xi^{1+\alpha} \left(\sin\left(\frac{\alpha \pi}{2}\right) - i \cos\left(\frac{\alpha \pi}{2}\right)  \right)} e^{i X \xi}  \, d\xi \\
%&+ \frac{1}{\sqrt{2 \pi}}   \int_{-\infty}^{0}  (-\xi)^s \, e^{-(-\xi)^{1+\alpha} \left(\sin\left(\frac{\alpha \pi}{2}\right) + i \cos\left(\frac{\alpha \pi}{2}\right)  \right)} e^{i X \xi} \, d\xi \\ 
=& \frac{1}{\sqrt{2 \pi}}   \int_{0}^{\infty}  \xi^\theta \, e^{-\xi^{1+\alpha} \left(\sin\left(\frac{\alpha \pi}{2}\right) - i \cos\left(\frac{\alpha \pi}{2}\right)  \right)} e^{- i (-X) \xi}  \, d\xi \\
&+ \frac{1}{\sqrt{2 \pi}}   \int_{0}^{\infty}  \xi^\theta \, e^{-\xi^{1+\alpha} \left(\sin\left(\frac{\alpha \pi}{2}\right) + i \cos\left(\frac{\alpha \pi}{2}\right)  \right)} e^{- i X \xi} \, d\xi,
\end{split}
\]
note that we have applied the change of variables $\xi \to -\xi$ in the second integral. Then, applying \cite[Lemma~2]{PT} in both integral terms, with $r=1+\alpha$, implies
\begin{equation}\label{D^sK:behaviour}
\left| |D|^\theta[K(1,\cdot)](X) \right| \lesssim \frac{1}{|X|^{1+\theta}}, \quad |X| \gg 1.
\end{equation}
Observe that the condition (\ref{cond:sigma}), since $\theta>0$, is satisfied for, in this case,
\[
\sigma= \sin\left(\frac{\alpha \pi}{2}\right) - i \cos\left(\frac{\alpha \pi}{2}\right).
\]
and its complex conjugate.

Then, (\ref{D^s:K}) and (\ref{D^sK:behaviour}), imply $|D|^\theta[K(1,\cdot)](X) \in L^p(\R)$ for $p\geq 1$, which together with (\ref{estimate:D^sK}) implies the third estimate.

Finally, the fourth estimate is obtained in a similar way. The main difference being that we have to differentiate the kernel first, which gives a factor $i\xi$ in the integrand, but we can still apply \cite[Lemma~2]{PT} to get
\[
%\begin{split}
| |D|^\theta [\partial_X K(1,\cdot)](X) | = \frac{1}{\sqrt{2\pi}}
\left| \int_{\R} |\xi|^\theta (i \xi) \, e^{(i\xi)^{1+\alpha}} e^{i X \xi} \, d\xi \right| 
% =& \frac{1}{\sqrt{2\pi}} \Bigg| i  \int_{0}^{\infty} \xi^{1+s} \, e^{-\xi^{1+\alpha} \left(\sin\left(\frac{\alpha \pi}{2}\right) - i \cos\left(\frac{\alpha \pi}{2}\right)  \right)} e^{i X \xi} \, d\xi \\  
% &-  i  \int_{-\infty}^{0}  (-\xi)^{1+s} \, e^{-(-\xi)^{1+\alpha} \left(\sin\left(\frac{\alpha \pi}{2}\right) + i \cos\left(\frac{\alpha \pi}{2}\right)  \right)} e^{i X \xi} \, d\xi \Bigg| \\
% \leq& \frac{1}{\sqrt{2\pi}} \left| \int_{0}^{\infty} \xi^{1+s} \, e^{-\xi^{1+\alpha} \left(\sin\left(\frac{\alpha \pi}{2}\right) - i \cos\left(\frac{\alpha \pi}{2}\right)  \right)} e^{-i (-X) \xi} \, d\xi \right| \\ 
% &+  \frac{1}{\sqrt{2\pi}}  \left|\int_{0}^{\infty}  \xi^{1+s} \, e^{-\xi^{1+\alpha} \left(\sin\left(\frac{\alpha \pi}{2}\right) + i \cos\left(\frac{\alpha \pi}{2}\right)  \right)} e^{- i X \xi} \, d\xi \right| \\
\lesssim \frac{1}{|X|^{2+\theta}}, \quad \text{for} \quad |X| \gg 1.
%\end{split}
\]
Then, with this and (\ref{K:self:sim}), we can argue as for (\ref{D^sK:estimate:t}) to conclude the proof.
\endproof

%---------------------------------------------------------------------------------------------------------------------------------------------------------------------

\subsection{Mild formulation, existence and regularity results}
%In this section, we recall some classical result for the purely convective problem (\ref{inviscid:Prob}) and for (\ref{NCL}). For the former we recall the results of \cite{LP}, first we give the definition of entropy solution in the sense of Kru\v{z}kov's (see \cite{Kruzhkov}):
%\begin{df}\label{entropy:SCL}
%Let $U_M$ be a weak solution of (\ref{inviscid:Prob}) such that
%\[
%U_M \in L^\infty((0,\infty), L^1(\R)) \cap L^\infty((\tau,\infty)\times \R), \ \forall \tau \in (0,\infty),
%\]
%then $U_M$ is said to be an entropy solution of (\ref{inviscid:Prob}) if, and only if, the following inequality holds for every $k\in \R$ and $\varphi \in C_c^\infty((0,\infty) \times \R)$ non-negative
%\begin{equation}\label{inviscid:entropy:inequality}
%\int_{0}^{\infty}\int_{\R} \left( |U_M - k| \partial_t \varphi + \sgn(U_M-k) (f(U_M) - f(k)) \partial_x \varphi \right) \, dx dt \geq 0,
%\end{equation}
%with $f(u) = |u|^{q-1} u/q$, and for any $\psi \in C_b(\R)$
%\begin{equation}\label{cons:mass:U:2}
%\limess_{t \downarrow 0} \int_{\R} U_M(t,x) \psi(x) \, dx = M \psi(0).
%\end{equation}
%\end{df}
%We recall that this unique entropy solution is given by the $N$-wave profile (see \cite{LP}).

%\begin{equation}\label{U:M}
%U_M(t,x) = 
%\begin{cases}
%\displaystyle\Big(\frac{x}{t}\Big)^{\frac{1}{q-1}} , \quad &0 < x < r(t), \\ \\
%\ \ \, 0, \quad &\text{otherwise},
%\end{cases}
%\ \ \mbox{with} \ \  r(t) = \left(\frac{q}{q-1}\right)^{\frac{q-1}{q}} \, M^{\frac{q-1}{q}} \, t^{\frac{1}{q}}.
%\end{equation}

We now define the mild formulation associated to (\ref{NCL}): 
%This mild formulation is just an equivalent formulation of (\ref{nonlinear:Prob}) which can be easily derived using the information for the linear fractional equation and Duhamel's principle. Therefore, one can arrive at the following definition of mild solution:
\begin{df}[Mild solution]\label{mild:solution}
Given $T\in (0,\infty]$ and $u_0 \in L^\infty(\R)$, we say that a mild solution of (\ref{NCL}) on $(0,T)\times\R$ is a function $u \in C_b((0,T)\times \R)$ which satisfies
\begin{equation}\label{mild:u}
u(t,x) = K(t,\cdot)*u_0(x) - \int_{0}^{t} \partial_x K(t-s,\cdot)*f(u(s,\cdot))(x) \, ds 
\end{equation}
in a.e. $(t,x) \in (0,T)\times \R$, where $f(u)= |u|^{q-1}u/q$ with $q>1$.
\end{df}

Regarding existence and uniqueness of mild solutions, we have the following:
\begin{thm}(Existence and uniqueness)\label{u:comparison:L1}
  Let $u_0\in L^\infty(\R)$, then there exists a unique global mild solution $u$ to the initial value problem (\ref{NCL}) with $u\in C((0, \infty),C^1(\R))\cap C_b((0, \infty)\times\R)$ and such that
\begin{equation}\label{max:min:prin}
  \essinf\{u_0\} \leq u(t,x) \leq \esssup\{u_0\}, \quad t>0\,,\quad x\in \mathbb{R}\,.
\end{equation}
If $u_0 \in L^\infty(\R) \cap L^1(\R)$, then also $u\in C([0,\infty), L^1(\R))$ and  %DO WE NEED THIS ONE?? DONT KNOW BUT MAYBE WE NEED THIS TO CONCLUDE THAT THE LIMIT SOLUTION U IS THE ENTROPY SOLUTION. I HAVE TO CHECK THIS
\begin{equation}\label{L1:con}
\|u(t,\cdot)\|_{L^1(\R)} \leq \|u_0\|_{L^1(\R)} \quad \mbox{for all} \quad t>0.
\end{equation}
\end{thm}

\proof
  The existence and uniqueness result, the upper bound of (\ref{max:min:prin}) and (\ref{L1:con}) have already been proved in \cite{DC} for a regular flux function.
  We observe that, in order to obtain existence and regularity we can proceed as in \cite[Propositions~4-5]{DC}, but since the flux function is only continuous with bounded first derivative, we can only apply two steps of the argument. This means that we can only gain $C^1$ regularity in $x$ and continuity for $t\in (0, T)$ for any $T>0$, thus $u\in C((0,T), C_b^1(\mathbb{R}))$.

 Now, in order to prove global existence, uniqueness and regularity as well as the upper bound analogous to the one in (\ref{max:min:prin}), we first regularise the flux function (to at least a $C^2$ function) and apply the results in the previous reference. The lower bound analogous to the one in (\ref{max:min:prin}) is proved following the proof of \cite[Lemma~3, Proposition~6]{DC} by changing the role of the supremum by the infimum, then we have to pass to the limit to get the results for the original flux. %We prove this in Appendix~\ref{app:1}.

We define the function $f_\delta$ by means of
\begin{equation}\label{reg:flux}
f_\delta(v) := \left( \delta^2 + v^2\right)^{\frac{q-1}{2}}  \frac{(v+\delta)}{q}.
\end{equation}
Notice that the function $f_\delta$ is $C^2$ for $\delta>0$ and converges to $f(v)=|v|^{q-1}v/q$ with $|f_\delta(v)-f(v)|\leq \delta C |v|^{q-1}$ as $\delta \to 0$, for some positive constant $C$.

Let $u_\delta$ be the solution of
\begin{equation}\label{NCL:nu}
\begin{cases}
\partial_t u_\delta(t,x) + \partial_x \left(f_\delta(u_\delta(t,x))\right) = \partial_x \DD^\alpha\left[u_\delta(t,\cdot)\right](x) ,& \quad t>0, \ x\in \R, \\
u_\delta(0,x) = u_0(x),& \quad x\in \R.
\end{cases}
\end{equation}
%We apply the results in \cite{DC} to conclude global existence, uniqueness and regularity as well as the upper bound analogous to the one in (\ref{max:min:prin}) for (\ref{NCL:nu}). The lower bound analogous to the one in (\ref{max:min:prin}) is proved following the proof of \cite[Lemma~3, Proposition~6]{DC} by changing the role of the supremum by the infimum. %We prove this in Appendix~\ref{app:1}.

As a consequence we obtain the global existence in time for (\ref{NCL:nu}). Therefore, by continuity in $t>0$ and the uniqueness result we can extend the solution for $t \in (0,\infty)$ and it satisfies
\[
\essinf {u_0} \leq u_\delta (t,x) \leq \esssup {u_0}
\]
for all $(t,x) \in (0,\infty)\times\R$. 

Now, we extend the result of (\ref{NCL:nu}) to (\ref{NCL}). First, we prove that for any $T>0$, $u_\delta$ converges uniformly to $u$ as $\delta \to 0$ in $(t,x) \in (0,T)\times\R$, where $u$ is a mild solution of (\ref{NCL}). We compute:
\[
\begin{split}
&\|u_\delta(t,\cdot) - u(t,\cdot)\|_{L^\infty(\R)} = \left\| \int_{0}^{t} \partial_x K(t-s,\cdot)*(f_\delta(u_\delta) - f(u))(x) \, ds \right\|_{L^\infty(\R)} \\
\leq& \int_{0}^{t} \|\partial_x K(t-s,\cdot)\|_{L^1(\R)} \| f_\delta(u_\delta(s,\cdot)) - f(u(s,\cdot)) \|_{L^\infty(\R)} \, ds \\
\leq& C \int_{0}^{t} (t-s)^{-\frac{1}{1+\alpha}}  \| f_\delta(u_\delta(s,\cdot)) - f(u_\delta(s,\cdot)) \|_{L^\infty(\R)} \, ds \\
&+ C \int_{0}^{t} (t-s)^{-\frac{1}{1+\alpha}}  \| f(u_\delta(s,\cdot)) - f(u(s,\cdot)) \|_{L^\infty(\R)} \, ds \\
%\leq&  \delta\, C(\|u_0\|_\infty) \int_{0}^{t} (t-s)^{-\frac{1}{1+\alpha}} \, ds \\
%&+ C L_f \int_{0}^{t} (t-s)^{-\frac{1}{1+\alpha}}  \| u_\delta(s,\cdot) - u(s,\cdot) \|_{L^\infty(\R)} \, ds \\
=&  \delta \, C(\|u_0\|_\infty)\, t^{\frac{\alpha}{1+\alpha}}  %\\ &
+ C(T) \int_{0}^{t} (t-s)^{-\frac{1}{1+\alpha}}  \| u_\delta(s,\cdot) - u(s,\cdot) \|_{L^\infty(\R)} \, ds,
\end{split}
\]
where the constant $C(T)$ depends linearly on $\sup_{t\in(0,T]}\|u(t,\cdot)\|_\infty$, which is finite for all $T$ (see \cite{DC}). Here we have used the second estimate with $p=1$ of Lemma~\ref{K:time:behaviour}, the convergence $f_\delta \to f$ as $\delta \to 0$ and the local Lipschitz continuity of $f$.

Since $1+\alpha>1$, we can apply the fractional Gronwall Lemma, see \cite[Lemma~2.4]{BC}, to obtain that for any $T>0$ and $\delta>0$ there exists a positive constant $C(T)$ such that
\[
\|u_\delta(t,\cdot) - u(t,\cdot)\|_{L^\infty(\R)} \leq \delta\, C(T), \quad \mbox{for all}\quad t \in [0,T].
\]
%This and the continuity of $u_\delta$ and $u$ with respect to $(t,x)$ imply the desired uniform convergence on $(0,T)\times\R$ for any $T>0$ fixed.
In particular, this implies that $u_\delta$ converges uniformly on compact sets of $t$ and point-wise (since $u\in C((0,T), C^1(\mathbb{R}))$. Thus, we can conclude global existence of $u$ and (\ref{max:min:prin}). 

Finally, (\ref{L1:con}) follows as in \cite[Theorem~2]{DC}.
\endproof

As a corollary we obtain positivity of the solutions for positive initial conditions:
\begin{cor}\label{Cor:pos_trans}
  Let $u_0\in L^\infty(\R)$ with $u_0(x) \geq \e >0$, then the unique mild solution of (\ref{NCL}) satisfies
\begin{equation}\label{positive}
\e \leq u(t,x) \leq \|u_0\|_\infty, \quad \mbox{for all} \quad t>0, \quad x\in \R.
\end{equation}
Moreover, $u\in C^\infty_b((0,\infty)\times \R)$.
\end{cor}
\proof
  The estimate (\ref{positive}) is a direct consequence of (\ref{max:min:prin}).

  Since now $u$ is positive this means that $|u|=u$ so the flux is $f(u) = u^q/q$ and belongs to $C^\infty([\e, \|u_0\|_\infty])$. This implies that then $u \in C^\infty_b((0,\infty)\times\R)$, see \cite{DC}.
\endproof

We now give some $L^p$-regularity of the mild solution.
%For more information on fractional Sobolev spaces see for example \cite{Adams}. The following regularity result is proved for the mild solution of (\ref{NCL}):
\begin{prop}[$L^p$-regularity]\label{u:p:regularity}
 Let $1< p < \infty$ and $u$ be the unique mild solution of (\ref{NCL}) with $u_0\in L^\infty(\R)\cap L^1(\R)$, then $\partial_t u \in C((0,\infty), L^p(\R))$ and $u \in C((0,\infty), L^p(\R) \cap\dot{H}^{\theta,p}(\R))$ for any $\theta< 1+\alpha + \min\{\alpha, q-1\}$.
\end{prop}
The proof is based on applying the fractional Laplacian to the mild formulation (\ref{mild:u}) followed by the bootstrap argument used in \cite[Proposition~3.1]{IS}. We use the fractional Laplacian to study regularity, because we can readily apply the chain and product rules available from \cite[Theorem~3, Corollary of Theorem~5]{Gatto}, \cite[Proposition 3.1]{ChW} and \cite[Proposition~A.1]{Visan}. This is also the reason why we need Lemma~\ref{K:time:behaviour} (the estimates needed to gain regularity by applying fractional derivatives on the kernel are obtained by convenience applying the fractional Laplacian instead of (\ref{WMFD})). With this ingredients one can just mimic the proof given in \cite{IS}.

Finally, we give another auxiliary result that we will need later on. Namely, the mild solution of (\ref{NCL}) satisfies a weak entropy inequality for the Kru\v{z}kov's entropies (see \cite[Theorem~1]{DC}):
\begin{thm}[Weak viscous entropy inequality]\label{weak:entropy:inequality}
For all $k\in \R$, let $\eta_k(v)=|v-k| \in C(\R)$ be a convex entropy function and $u \in C((0,\infty),C^1(\R))\cap C_b((0,\infty)\times \R)$ a solution of (\ref{NCL}), then for all non-negative $\varphi \in C_c^\infty ((0,\infty)\times \R)$
\begin{equation}
\begin{split}
\int_{0}^{\infty} \int_{\R} \Big( |u(t,x) - k| \partial_t \varphi + \sgn(u(t,x)-k)(f(u(t,x))-f(k)) \partial_x \varphi& \\
+ |u(t,x) - k| \partial_x \overline{\DD^\alpha}[\varphi(t,\cdot)](x)  \Big) \, dx dt  &\geq 0,
\end{split}
\label{Kruzhkov:entropy:inequality}
\end{equation}
where $f(u)=|u|^{q-1}u/q$.
\end{thm} 
We remark that the proof of this result is as in \cite{DC}, but for the problem with the regularised flux (\ref{reg:flux}), as above. The proof is completed by passing to the limit $\delta\to 0$. We omit the details here.

%---------------------------------------------------------------------------------------------------------------------------------------------------------------------

\section{Oleinik type inequality for non-negative solutions}\label{sec:3}
In this section we derive 
%As mentioned before, the key estimate in the proof of asymptotic behaviour of solution is
an Oleinik type inequality. We prove it for non-negative solutions by first deriving the inequality for positive ones (for which the flux is regular).
%Thus, we approximate the solution by a family of positive solutions of an approximating problem. As a consequence of the previous Corollary~\ref{Cor:pos_trans} the approximated solution will preserve positivity and after proving the necessary estimates we will pass to the limit and obtain the inequalities for the initial solution.

Let $u_0\in L^\infty(\R)$ be non-negative, then we consider the following approximating problem,
\begin{equation}\label{approx:Prob}
\begin{cases}
\partial_t u_\e(t,x) +  (u_\e)^{q-1} \partial_x u_\e(t,x) = \partial_x \DD^\alpha[u_\e(t,\cdot)](x), & \quad t>0,  \ x\in\R, \\
u_\e(0,x) = u_0(x) + \e, & \quad x\in\R.
\end{cases}
\end{equation}
The existence and uniqueness of this problem is guaranteed by Theorem~\ref{u:comparison:L1} and Corollary~\ref{Cor:pos_trans} implies that $u_\e(x,t)>0$ for all $(x,t)$. Then, the following holds:
\begin{lem}\label{eps:limit}
 Let $u(t,x)$ be the solution of (\ref{NCL}) with $0\leq u_0 \in L^\infty(\R)$ and let $u_\e(t,x)$ be the solution of (\ref{approx:Prob}). Then for every $T>0$,
\begin{equation*}
\max_{t\in [0,T]} \|u_\e(t,\cdot) - u(t,\cdot)\|_{L^\infty(\R)} \to 0 \quad \text{as} \quad \e\to 0. 
\end{equation*} 
\end{lem}
The proof is analogous to that in \cite{IS}, we do not prove it here. It is based on the comparison of the mild solutions in norm and on the application of the fractional Gronwall Lemma \cite{BC}.

%---------------------------------------------------------------------------------------------------------------------------------------------------------------------
We can now obtain an Oleinik type entropy inequality. 
\begin{prop}[Oleinik entropy inequality]\label{Oleinik:inequality}
Let $u_0\in L^\infty(\R)$. Then, for any $\e>0$, the solution $u_\e$ of (\ref{approx:Prob}) satisfies
\begin{equation*}
\partial_x \left(u_\e^{q-1}(t,x)\right) \leq \frac{1}{t}, \quad \forall t>0, \ x\in \R. 
\end{equation*}
\end{prop}

The proof is analogous to that in \cite{IS}, we briefly give the idea here since this proposition is key to the proof of Theorem~\ref{subcritical:asymptotics}: First, one writes the equation for $w=(u_\e^{q-1})_x$, 
\begin{equation}\label{w:eq}
w_t + w^2 + u_\e^{q-1} \, w_x + N(w,u_e)=0,
\end{equation}
where $N$ is the resulting non-linear and non-local term,
and shows that $W=\sup_{x\in\R} {w(\cdot,x)}$, satisfies
\[
W'(t) + W^2(t) +Q(t) W(t)\leq 0, \quad Q(t)>0
\]
which, by ODEs arguments, implies the result. The estimates of the non-local term  of (\ref{w:eq}) can be adapted from the proof in \cite{IS}. This requires splitting the integral, using the representation (\ref{int:frac:Lap}), near zero and away from zero, but this can be done similarly using the integral representations of Lemma~\ref{equiv:repre} in our case.

%---------------------------------------------------------------------------------------------------------------------------------------------------------------------

We can show properties of the family $u_\e$ to the solution of (\ref{NCL}) by taking the limit $\e\to 0$. Namely,
\begin{lem}\label{properties:u} Let $u$ be the solution of problem (\ref{NCL}) with non-negative initial data $u_0 \in L^\infty(\R)\cap L^1(\R)$. Then, the following estimates hold:
\begin{enumerate}
\item[(i)] (Mass conservation) $\int_{\R}{u(t,x) \, dx}=M$, $\forall t>0$ where $M$ is defined as $M= \int_{\R}{u_0(x)\, dx}$.
\item[(ii)] (Oleinik entropy condition) $\partial_x\left( u^{q-1}(t,x)\right) \leq 1/t$ for all $t>0$ in a weak distributional sense.
\item[(iii)] (Upper bound) $0 \leq u(t,x) \leq \left( \frac{q}{q-1} M\right)^{1/q} \, t^{-1/q}$ for all $t>0$ and $x\in \R$.
\item[(iv)] (Decay in $L^p$-norm) For $1\leq p \leq \infty$,
\[
\| u(t,\cdot)\|_{L^p(\R)} \leq \left(\frac{q}{q-1}\right)^{\frac{p-1}{pq}} \, M^{\frac{p-1}{pq} + \frac{1}{p}} \, t^{-\frac{1}{q} (1- \frac{1}{p})}, \quad \forall t>0.
\]
\item[(v)] (Decay of the spatial derivative) $\partial_x u(t,x) \leq C(q) M^{\frac{2-q}{q}} \, t^{-\frac{2}{q}}$, for all $t>0$ and a.e. $x\in \R$.

\item[(vi)] ($W^{1,1}_{loc}(\R)$ estimate) For any $R>0$,
\[
\int_{|x|<R}{\left| \partial_x u(t,x)\right| \, dx} \leq 2R\, C(q)\, M^{\frac{2-q}{q}} \, t^{-\frac{2}{q}} + 2\left(\frac{q}{q-1} M\right)^{\frac{1}{q}} \, t^{-\frac{1}{q}}, \quad \forall t>0.
\]

\item[(vii)] (Energy estimate) For every $0<\tau < T$,
\[
\int_{\tau}^{T} \int_{\R}{\left| \DD^{\frac{\alpha+1}{2}}[u(t,\cdot)](x)\right|^2 \, dx dt}
%\leq \frac{1}{2} \int_{\R}{u^2(\tau , x) \, dx}
\leq
\frac{1}{2} \left(\frac{q}{q-1}\right)^{\frac{1}{q}} \, \tau^{-\frac{1}{q}} \, M^{\frac{q+1}{q}}.
\]
\end{enumerate}
\end{lem}

\proof
The proof of (i) is as in \cite{IS}.

We recall the proof of (ii):  First we recall that $u \in C((0,\infty), L^p(\R) \cap\dot{H}^{\theta,p}(\R))$ for any $\theta<1 + \alpha+\min\{\alpha,q-1\}$ and $u_\e \in C_b^\infty((0,\infty)\times \R)$, then using Lemma~\ref{eps:limit}, $u_\e \to u$ as $\e\to 0$ point-wise for all $t>0$ and $x\in \R$. As a result, (ii)  holds by taking the limit $\e\to 0$ and Proposition~\ref{Oleinik:inequality}. With this, one can again proceed as in \cite{IS} to conclude (iii) and (iv). The proofs of (v) and (vi) follow as in \cite{IS} too, they do not depend on the form of the non-local operator.

We now prove (vii).  First, we multiply (\ref{NCL}) by $u$ and get the following identity, after integrating with respect to $x$,
\[
%\label{energy:estimate:1}
\frac{1}{2}\frac{d}{dt} \int_{\R} u^2 \, dx - \int_{\R}u \partial_x \DD^\alpha[u](x) \, dx = 0.
\]
We observe that, using integration by parts
\[
- \int_{\R}u \partial_x \DD^\alpha[u](x) \, dx = \int_{\R} \partial_xu\,  \DD^\alpha[u](x) \, dx \geq 0,
\]
the last inequality is shown in \cite[Lemma~2.2]{CA}. Then, proceeding as in the proof of Lemma~\ref{partial:int:parts:D:alpha} with $\theta_1=\theta_2$, using (\ref{f:symbol:D:alpha})-(\ref{f:symbol:adjoint:D:alpha}) with $\alpha$ replaced by $(\alpha+1)/2$, we have
\[
\begin{split}
0\leq -\int_{\R}u \partial_x \DD^\alpha[u](x) \, dx = \int_{\R} (i \xi) \FF(u)(\xi) \overline{(i \xi)^{\alpha} \FF(u)(\xi)} \, d\xi\\
=\left| \int_{\R} (-1)^{\frac{1+\alpha}{2}+1} \, (i \xi)^{\frac{1+\alpha}{2}}  (-i \xi)^{\frac{1+\alpha}{2}} \FF(u)(\xi) \overline{\FF(u)(\xi)} \, d\xi \right| \\
=  \int_{\R} (i \xi)^{\frac{1+\alpha}{2}}\FF(u)(\xi) \overline{(i \xi)^{\frac{1+\alpha}{2}}\FF(u)(\xi)} \, d\xi. 
\end{split}
\]
%the imaginary part of the first integral is zero using the symmetry
%\[
%|\FF(u)(\xi)|^2 = |\FF(u)(-\xi)|^2.
%\]
Applying again Plancherel's theorem we finally have
\begin{equation}\label{energy:estimate:2}
\frac{1}{2} \frac{d}{dt}  \int_{\R} u^2(t,x) \, dx  + \int_{\R}{\left| \DD^{\frac{\alpha+1}{2}}[u(t,\cdot)](x)\right|^2 \, dx}  = 0.
\end{equation}
The conditions to be able to apply Plancherel's theorem follow from (iv) and Proposition~\ref{u:p:regularity}.

Finally, we integrate (\ref{energy:estimate:2}) over $(\tau,T)$ for some $\tau>0$, to get
\[
\frac{1}{2}  \int_{\R} u^2(T,x) \, dx - \frac{1}{2}  \int_{\R} u^2(\tau,x) \, dx + \int_{\tau}^{T}\int_{\R}{\left| \DD^{\frac{\alpha+1}{2}}[u(t,\cdot)](x)\right|^2 \, dx \, dt}  = 0.
\]
Then, (vii) follows taking into account that the first term is non-negative and applying (iv) for $p=2$.
\endproof

%---------------------------------------------------------------------------------------------------------------------------------------------------------------------

\section{Asymptotic behaviour}\label{sec:4}%{asymptotic:behaviour}
In this section we prove Theorem~\ref{subcritical:asymptotics}. As we have mentioned, this is equivalent to proving the limit (\ref{asymp:ulambda}) for some fixed $s$. Henceforth, we consider for $\lambda>1$, $u_\lambda$ be defined by means of (\ref{lambda:scaling})-(\ref{u:scaling}).% is d
We follow the proof in \cite{IS}, we use the Lemma~\ref{properties:u} (translated to $u_\lambda$) and the following lemma:

\begin{lem}[Tail control estimate]\label{properties:u:lambda}
Let $\lambda>1$ and $u_\lambda$ be the solution of (\ref{lambda:Prob}) with $0\leq u_0 \in L^\infty(\R) \cap L^1(\R)$. Then,
for any $\lambda>1$ and $R>0$, there exists a constant $C(M,q)>0$ such that
\[
\int_{|y|>2R} u_\lambda(s,y) \, dy \leq \int_{|y|>R} u_0(y) \, dy + C(M,q) \left( \frac{s \lambda^{q-1-\alpha}}{R^{\alpha+1}} + \frac{s^{1/q}}{R} \right), \quad \forall s>0.
\]
\end{lem} 
%WHERE DO WE USE THIS LEMMA?%\proof
Noticing that $\DD^{\alpha}[u_\lambda(s, \cdot)](y) = \lambda^{\alpha+1} \, \DD^{\alpha}[u(\lambda^q s, \cdot)](\lambda y)$ and the same scaling works for $\partial_y \overline{\DD^{\alpha}}[\cdot]$, we can just follow \cite{IS} to conclude the result.

\proof[Proof of Theorem~\ref{subcritical:asymptotics}]
The existence and the regularity item (i) and (ii) follow from Theorem~\ref{u:comparison:L1} and Proposition~\ref{u:p:regularity}. It remains to show the large time behaviour.

First, we show that (\ref{asymp:sub}) and (\ref{asymp:ulambda}) are equivalent. Without loss of generality we consider $s_0=1$, applying the scaling (\ref{lambda:scaling})-(\ref{u:scaling}), and observing that (\ref{inviscid:Prob}) is invariant under the rescaling, we get that for any $\lambda>1$,
\[
%\label{ulambda:u}
u_\lambda(1,y) - U_M(1,y) =t^{\frac{1}{q}} \left(u(t , x) - U_M(t,x)\right).
\]
And performing the change of variables in the integral, we finally get
\[
\|u_\lambda(1,y) - U_M(1,y)\|_{L^p(\R)} = t^{\frac{1}{q}\left(1-\frac{1}{p}\right)} \, \|u(t,x) - U_M(t,x)\|_{L^p(\R)}.
\]

Let us then prove (\ref{asymp:ulambda}). We divide the proof into several steps. Let us first show the convergence of a sub-sequence of $\{u_\lambda\}_{\lambda>1}$. Using, \cite[Theorem~5]{Simon} we shall get that $\{u_\lambda\}_{\lambda>1}$ is relatively compact in $C([s_1,s_2], L^2_{loc}(\R))$ for any $0< s_1 < s_2 < \infty$.

We let $B_R = (-R,R)$ and we apply \cite[Theorem~5]{Simon} to the triple $W^{1,1}(B_R) \hookrightarrow L^2(B_R) \hookrightarrow H^{-1}(B_R)$. Observe that (i) and (iii) in Lemma~\ref{properties:u}, imply that $\{u_\lambda\}_{\lambda>1}$ is uniformly bounded in $L^\infty((s_1,s_2), W^{1,1}(B_R))$, and this gives the first condition of this theorem. Then, by \cite[Lemma~4]{Simon} we can conclude that
\begin{equation}\label{2:simon}
\| u_\lambda(s+h, \cdot) - u_\lambda(s, \cdot) \|_{L^\infty((0,T-h), H^{-1}(B_R))} \to 0 \quad \text{as} \quad h\to 0 \ \text{uniformly for} \ \lambda>1
\end{equation}
provided that $\{ \partial_s u_\lambda \}_{\lambda>1}$ is uniformly bounded in $L^p((s_1,s_2), H^{-1}(B_R))$ for some $p<\infty$. 
Let us show this with $p=2$. First, let us choose $\varphi \in C_c((0,\infty)\times B_R)$ and extend it by zero outside $B_R$, for such $\varphi$ and $\lambda > 1$ we have
\begin{equation}\label{u:lambda:weak:ineq}
\begin{split}
\Bigg| \int_{s_1}^{s_2}  \int_{\R} &(\partial_s u_\lambda)  \, \varphi \, dy ds \Bigg| \leq  \left| \int_{s_1}^{s_2}  \int_{\R} \partial_y (u_\lambda)^q \, \varphi \, dy ds \right|  + \lambda^{q-1-\alpha}  \left| \int_{s_1}^{s_2}  \int_{\R} \partial_y \DD^\alpha[u_\lambda] \, \varphi \, dy ds \right| \\
=& \left| \int_{s_1}^{s_2}  \int_{\R} u_\lambda^q \, \partial_y \varphi \, dy ds \right|  +  \lambda^{q-1-\alpha} \left| \int_{s_1}^{s_2} \int_{\R} \DD^{\frac{1+\alpha}{2}}[u_\lambda] \, \overline{\DD^{\frac{1+\alpha}{2}}}[\varphi] \, dy ds \right|  \\
\leq& \left\| u_\lambda^q \right\|_{L^2((s_1,s_2), L^2(\R))} \, \| \varphi \|_{L^2((s_1,s_2), H^1(\R))}  \\
 &+   \lambda^{\frac{q-1-\alpha}{2}} \left( \lambda^{q-1-\alpha} \int_{s_1}^{s_2} \int_{\R} \left| \DD^{\frac{1+\alpha}{2}}[u_\lambda] \right|^2 \, dy ds\right)^{\frac{1}{2}} \left( \int_{s_1}^{s_2} \int_{\R} \left| \overline{\DD^{\frac{1+\alpha}{2}}}[\varphi] \right|^2 \, dy ds \right)^{\frac{1}{2}} \\
 \leq& C'(M,q,s_1,s_2)  \| \varphi \|_{L^2((s_1,s_2), H^1(\R))}  \\
 &+   \lambda^{\frac{q-1-\alpha}{2}} \frac{1}{\sqrt{2}} \left(\frac{q}{q-1}\right)^{\frac{1}{2q}} M^{\frac{1+q}{2q}} s_1^{-\frac{1}{2q}} \| \varphi\|_{L^2\left((s_1,s_2), \dot{H}^{\frac{1+\alpha}{2}}(\R)\right)} \\
 \leq& C(M,q,s_1,s_2)  \| \varphi \|_{L^2((s_1,s_2), H^1(\R))}.
\end{split}
\end{equation}
Here, we have applied Lemma~\ref{partial:int:parts:D:alpha} in the second inequality and the energy estimate Lemma~\ref{properties:u} (iv), as well as (\ref{homo:norm}). All these steps can be performed since conservation of mass and the regularity of $u$ is transferred to $u_\lambda$ (see Proposition~\ref{u:p:regularity}, in particular) and by the choice of $\varphi$. Now, the Riesz representation  theorem and \cite[Chapter~IV Corollary~4]{DU}, (\ref{u:lambda:weak:ineq}) imply that
\[
\| \partial_s u_\lambda \|_{L^2((s_1,s_2), H^{-1}(B_R))}  \leq C(M,q,s_1,s_2), \quad \forall \lambda > 1,
\]
and we can conclude (\ref{2:simon}). Hence, we can apply \cite[Theorem~5]{Simon}, this means that $\{ u_\lambda \}_{\lambda>1}$ is relatively compact in $C([s_1,s_2], L^2(B_R))$.

As a consequence there exists $u_\infty \in C([s_1,s_2], L^2(B_R))$ such that, up to a sub-sequence, $u_\lambda \to u_\infty$ as $\lambda \to \infty$ in  $C([s_1,s_2], L^2(B_R))$.
By a diagonal argument we can conclude the convergence for any compact set and, therefore, 
\begin{equation}\label{local:L2:convergence}
u_\lambda \longrightarrow u_\infty, \quad \text{as} \ \ \lambda \to \infty \ \ \text{in} \ \ C([s_1,s_2], L^2_{loc}(\R)).
\end{equation}
We observe that (\ref{local:L2:convergence}) implies also that $u_\lambda \to u_\infty$ in $C([s_1,s_2], L^1_{loc}(\R))$. In order to extend this convergence to $C([s_1,s_2], L^1(\R))$, we use Lemma~\ref{properties:u:lambda} (see \cite{IS}). Hence,
\[
u_\lambda \longrightarrow u_\infty, \quad \text{as} \ \ \lambda \to \infty \ \ \text{in} \ \ C([s_1,s_2], L^1(\R)).
\]

The next step is to prove that $u_\infty=U_M$, i.e. that it satisfies Definition~\ref{entropy:SCL}. First, we recall that $u$ satisfies (\ref{Kruzhkov:entropy:inequality}) of Theorem~\ref{weak:entropy:inequality}. Therefore, $u_\lambda$ satisfies the following inequality for any non-negative $\varphi \in C_c^\infty ((0,\infty)\times \R)$, using Lemma~\ref{equiv:repre}:
\begin{equation}\label{wvei}
\begin{split}
\int_{0}^{\infty} \int_{\R} \Big( |u_\lambda - k| \partial_s \varphi & +\frac{1}{q} \sgn(u_\lambda - k) \, ((u_\lambda)^q - k^q) \partial_y \varphi  \\
&+ \lambda^{q-1-\alpha} \, \partial_y \DD^\alpha [|u_\lambda - k|](y) \varphi   \Big) \, dy ds \geq 0.
\end{split}
\end{equation}
In what follows we pass to the limit $\lambda \to\infty$ in (\ref{wvei}).
We prove that the last term tends to zero as $\lambda \to \infty$. We split this integral term into two as follows, given $r>0$,
\[
\begin{split}
 &\int_{0}^{\infty}\!\! \int_{\R}\partial_y \DD^\alpha [|u_\lambda - k|](y) \varphi(s,y)dy ds \\& =d_{\alpha+2} \int_{0}^{\infty} \!\!\int_{\R}  \!\! 
  \int_{-\infty}^{-r} \frac{ |u_\lambda(s,y+z) - k| - |u_\lambda(s,y) - k| - \partial_y (|u_\lambda - k|) z }{|z|^{\alpha+2}}\varphi(s,y) \, dz   dy ds \\
 & +d_{\alpha+2} \int_{0}^{\infty} \!\!\int_{\R}  | u_\lambda(s,y) - k| \int_{0}^{r} \frac{ \varphi(s,y+z) - \varphi(s,y) - \partial_y \varphi \, z }{|z|^{\alpha+2}} \, dz dy ds \, . % \\
% &= I_1 + I_2.  
\end{split}
\] 
The second integral term has been obtained by using Fubini's theorem, integration by parts in $y$ in the third term, and the pertinent changes of variables.

Following the ideas of \cite{IS}, we bound the first and second integral terms applying the regularity of $\varphi$ and the non-negativity and conservation of mass of $u_\lambda$. Then, this means that the last term in the left-hand side of (\ref{wvei}) goes to zero as $\lambda\to\infty$.

Since $u_\lambda \to u_\infty$ in $C((0,\infty), L^1(\R))$, we can pass to the limit in property (i) of Lemma~\ref{properties:u}, so that $\int_{\R} u_\infty(s,y) \, dy = M$. Moreover, $u_\lambda \to u_\infty$ a.e. in $(0,\infty)\times \R$ which shows that property (ii) of Lemma~\ref{properties:u} with $p=\infty$ is transferred to $u_\infty$:
\[
\|u_\infty(s, \cdot)\|_{L^\infty(\R)} \leq C(M) s^{-\frac{1}{q}}.
\]
This last inequality is sufficient to prove that $(u_\lambda)^q \to (u_\infty)^q$ as $\lambda \to \infty$ in $C((0,\infty), L^1(\R))$ and, therefore, passage to the limit $\lambda \to\infty$ in (\ref{wvei}) gives Definition~\ref{entropy:SCL}-(\ref{inviscid:entropy:inequality}) (with $U_M$ replaced by $u_\infty$) for every constant $k\in \R$ and $\varphi \in C_c^\infty ((0,\infty)\times \R)$, $\varphi \geq 0$. 

Finally, we have to check that $u_\infty$ satisfies Definition~\ref{entropy:SCL}-(\ref{cons:mass:U:2}) for any $\psi\in C_b(\R)$. First, one proves it for any $\psi \in C^2_b(\R)$, which, by density, can be generalised to $\psi\in H^2(\R)$. Finally, the result with $\psi \in C_b(\R)$ follows by an approximation argument and Lemma~\ref{properties:u}.

 Thus we have shown that $u_\infty$ satisfies Definition~\ref{entropy:SCL}. Since (\ref{inviscid:Prob}) has a unique entropy solution,  $U_M$, then $\{ u_\lambda \}_{\lambda>1}$ converges to $U_M$ in $C([s_1,s_2], L^1(\R))$ as $\lambda\to\infty$.

% we have proved for $p=1$ that $u_\lambda(s,\cdot) \to U_M(s,\cdot)$ as $\lambda \to \infty$ in $L^1(\R)$ for any $s>0$.
 In order to finish the proof, one has to extend this convergence to $L^p(\R)$ with $1 < p < \infty$. This follows by interpolation as in \cite{IS}.

 %Indeed, applying interpolation for $\frac{1}{p} = \frac{1-\theta}{1} + \frac{\theta}{2p}$ with $\theta = \frac{2(p-1)}{2p-1}$, we get
%\[
%\begin{split}
%\| u_\lambda(s,\cdot) - U_M(s,\cdot) \|_{L^p(\R)} \leq& \| u_\lambda(s,\cdot) - U_M(s,\cdot) \|_{L^1(\R)}^{1/(2p-1)} \\
%&\cdot \left( \| u_\lambda(s,\cdot) \|_{L^{2p}(\R)}  + \| U_M(s,\cdot) \|_{L^{2p}(\R)}\right)^{2(p-1)/(2p-1)}.
%\end{split}
%\]
%And passing to the limit $\lambda\to\infty$ we obtain the result.
\endproof

%----------------------------------------------------------------------------------------------------------------------------------------------------------------------------------------------------------------

\section{Regularisation by a general Riesz-Feller operator}\label{sec:5}%{asymp:general:RF}
In this section, we focus on showing how to generalise the previous results of Sects.~\ref{sec:2}, \ref{sec:3} and \ref{sec:4} for the problem
\begin{equation}\label{RFCL}
\begin{cases}
\partial_t u(t,x) + |u(t,x)|^{q-1}\partial_x (u(t,x)) = D^\beta_\gamma[u(t,\cdot)](x), \quad &t>0, \ \ x\in \R, \\
u(0,x) = u_0(x), \quad &x\in \R,
\end{cases}
\end{equation}
where the diffusion is given by a general Riesz-Feller operator (\ref{Riesz-Feller:operator})-(\ref{Riesz-Feller:symbol}). Here $\beta$ and $\gamma$ satisfy the assumptions of such definition and $u_0 \in L^\infty(\R)\cap L^1(\R)$.

We use the following formulation of the non-local operator, given in \cite[Proposition~2.3]{AK2} (or see \cite{Cifani,MLP,Sato}): for any $0<\beta<2$ and $|\gamma| \leq \min\{ \beta, 2 - \beta\}$,
\[
\begin{split}
D^\beta_\gamma[g](x) =& c^1_\gamma \int_{0}^{\infty} \frac{g(x-z) -g(x) + g'(x) z}{z^{1+\beta}} \, dz \\
 &+ c^2_\gamma \int_{0}^{\infty} \frac{g(x+z) -g(x) - g'(x) z}{z^{1+\beta}} \, dz, \quad \text{for} \ 1<\beta<2,
\end{split}
%\label{D:beta:1:equiv}
\]
where (e.g. see \cite{MLP}) 
\[
c^1_\gamma = \frac{\Gamma(1+\beta)}{\pi}\sin\left( (\beta-\gamma) \frac{\pi}{2}\right) \quad  \text{and} \quad c^2_\gamma = \frac{\Gamma(1+\beta)}{\pi}\sin\left( (\beta+\gamma) \frac{\pi}{2}\right),
\]
in particular $c^1_\gamma + c^2_\gamma >0$.

Using Lemma~\ref{equiv:repre} it is easy to show that%one can show that
\begin{equation}\label{RF:integral:representation}
D^\beta_\gamma[ g](x)= \frac{1}{d_{\beta+1}} \left(c_\gamma^1\partial_x\DD^{\beta-1}[g](x)
+ c_\gamma^2 \partial_x\overline{\DD^{\beta-1}}[g](x)\right).
\end{equation}
%for $c^1_\gamma, c^2_\gamma \geq 0$ with $c^1_\gamma + c^2_\gamma > 0$.

Existence and regularity results for (\ref{RFCL}) are proved similarly by defining mild solutions as in Definition~\ref{mild:solution} with the kernel
\begin{equation}
\label{kernel:RFP}
K^\beta_\gamma(t,x) := \FF^{-1}\left( e^{t \psi^\beta_\gamma(\xi)} \right)(x),
\end{equation}
these steps have already been explained in \cite[\S~6]{DC}. This is because, $K^\beta_\gamma$ satisfies similar properties as $K$ does, the proofs are given in e.g. \cite[Lemma~2.1]{AK2}. Thus, we can say that Theorem~\ref{u:comparison:L1} and Corollary~\ref{Cor:pos_trans} hold unchanged for (\ref{RFCL}).

In order to generalise Proposition~\ref{u:p:regularity} for (\ref{RFCL}), we need the following lemma (analogous to Lemma~\ref{K:time:behaviour}):

\begin{lem}[Time behaviour of $K^\beta_\gamma$]\label{K:RF:time:behaviour} 
For all $\theta \in (0,1)$ and $1\leq p \leq \infty$, the kernel $K^\beta_\gamma(t,x)$, such that $\beta \in (1,2)$ and $|\gamma| \leq \min\{ \beta, 2-\beta\}$, satisfies the following estimates for any $t>0$:
\[
\begin{split}
\|K^\beta_\gamma(t,\cdot)\|_{L^p(\R)}& = C t^{-\frac{1}{\beta} (1-\frac{1}{p})}, \\ 
\|\partial_x K^\beta_\gamma(t,\cdot)\|_{L^p(\R)}& \lesssim t^{-\frac{1}{\beta} (1-\frac{1}{p}) - \frac{1}{\beta}}, \\
\| |D|^\theta [K^\beta_\gamma(t,\cdot)]\|_{L^p(\R)}& \lesssim  t^{-\frac{1}{\beta} (1-\frac{1}{p})- \frac{\theta}{\beta}}, \\ 
\| |D|^\theta [\partial_x K^\beta_\gamma(t,\cdot)]\|_{L^p(\R)}& \lesssim  t^{-\frac{1}{\beta} (1-\frac{1}{p})- \frac{1+\theta}{\beta}}, \\ 
\end{split}
\] 
for some constant $C>0$.
\end{lem}
\proof
The properties of (\ref{kernel:RFP}), analogous to the ones for (\ref{kernel:CP}) are given in \cite[Lemma~2.1]{AK2}, and combining the self-similarity, the mass conservation of (\ref{kernel:RFP}) and its derivative, and the fact that these are bounded on $(0,T)\times\R$ for any $T>0$, we conclude the first and second estimates.

For the third estimate, we apply the self-similarity property of $K^\beta_\gamma$ and rescale as follows% scaling property of the fractional Laplacian:
\begin{equation}\label{D^sK:RF:estimate:t}
| |D|^\theta[K^\beta_\gamma(t,\cdot)](x)| = \frac{1}{t^{1/\beta}} \left| |D|^\theta\left[K^\beta_\gamma\left(1,\frac{\cdot}{t^{1/\beta}}\right)\right](x)\right| = \frac{1}{t^{\frac{1+\theta}{\beta}}} \left| |D|^\theta[K^\beta_\gamma(1,\cdot)]\left( \frac{x}{t^{1/\beta}}\right)\right|.
\end{equation}
If the $L^p$-norm of (\ref{D^sK:RF:estimate:t}) is finite, we get, applying the change of variable $X= \frac{x}{t^{1/\beta}}$, the desired estimate:
\[
%\label{estimate:D^sK:RF}
\begin{split}
  \|  |D|^\theta[K^\beta_\gamma(t,\cdot)](x) \|_{L^p(\R)}
  %= \frac{1}{t^{\frac{1+\theta}{\beta}}} \left(\int_{\R} \left| |D|^\theta[K^\beta_\gamma(1,\cdot )]\left(\frac{x}{t^{\frac{1}{\beta}}}\right)\right|^p \, dx \right)^{1/p} 
= \frac{1}{t^{\frac{1+\theta}{\beta}}} t^{\frac{1/p}{\beta}}  \left(\int_{\R} \left| |D|^\theta[K^\beta_\gamma(1,\cdot )](X)\right|^p \, dX \right)^{1/p}  
\lesssim  t^{-\frac{1}{\beta} \left(1-\frac{1}{p}\right)- \frac{\theta}{\beta}}.
\end{split}
\]

Thus, it remains to show that the $L^p$-norm is finite. One gets the boundedness of the integrand using the definition (\ref{Fractional:Laplacian})
\begin{equation}\label{D^s:K:RF}
\begin{split}
  \left| |D|^\theta[K^\beta_\gamma(1,\cdot )] (X)\right|
  =& \frac{1}{\sqrt{2\pi}} \left| \int_{\R} |\xi|^\theta \, e^{-|\xi|^\beta e^{i \sgn(\xi) \gamma \frac{\pi}{2}}} e^{i X \xi} d\xi \right| \\
 \leq& \frac{1}{\sqrt{2\pi}} \int_{\R} |\xi|^\theta \, e^{-|\xi|^{\beta} \cos\left( \frac{\gamma \pi}{2}\right)} \, d\xi < \infty,
\end{split}
\end{equation}
where $|\gamma| \leq 2-\beta < 1$, which implies that $\cos\left( \frac{\gamma \pi}{2}\right)>0$. Hence, in order to conclude, it is sufficient to control the behaviour for large $|X|$. Starting from (\ref{D^s:K:RF}), we write,
\[
\begin{split}
|D|^\theta[K^\beta_\gamma(1,\cdot)](X) =& \frac{1}{\sqrt{2 \pi}}  \int_{\R} |\xi|^\theta \, e^{-|\xi|^\beta e^{i \sgn(\xi) \gamma \frac{\pi}{2}}}  e^{i X \xi} \, d\xi \\
=&  \frac{1}{\sqrt{2 \pi}}  \int_{\R} |\xi|^\theta \, e^{-|\xi|^{\beta} \left(\cos\left(\frac{\gamma \pi}{2}\right) + i \sgn(\xi)\sin\left(\frac{\gamma \pi}{2}\right)  \right)} e^{i X \xi} \, d\xi \\
%=& \frac{1}{\sqrt{2 \pi}}  \int_{0}^{\infty}  \xi^s \, e^{-\xi^{\beta} \left(\cos\left(\frac{\gamma \pi}{2}\right) + i \sin\left(\frac{\gamma \pi}{2}\right)  \right)} e^{i X \xi}  \, d\xi \\
%&+ \frac{1}{\sqrt{2 \pi}}   \int_{-\infty}^{0}  (-\xi)^s \, e^{-(-\xi)^{\beta} \left(\cos\left(\frac{\gamma \pi}{2}\right) - i \sin\left(\frac{\gamma \pi}{2}\right)  \right)} e^{i X \xi} \, d\xi \\ 
=& \frac{1}{\sqrt{2 \pi}}   \int_{0}^{\infty}  \xi^\theta \, e^{-\xi^{\beta} \left(\cos\left(\frac{\gamma \pi}{2}\right) + i \sin\left(\frac{\gamma \pi}{2}\right)  \right)} e^{- i (-X) \xi}  \, d\xi \\
&+ \frac{1}{\sqrt{2 \pi}}   \int_{0}^{\infty}  \xi^\theta \, e^{-\xi^{\beta} \left(\cos\left(\frac{\gamma \pi}{2}\right) - i \sin\left(\frac{\gamma \pi}{2}\right)  \right)} e^{- i X \xi} \, d\xi.
\end{split}
\]
%note that in the second integral the change of variable $\xi \to -\xi$ is applied.
Now, let 
\[
\sigma= \cos\left(\frac{\gamma \pi}{2}\right) + i \sin\left(\frac{\gamma \pi}{2}\right).
\] 
Now we apply \cite[Lemma~2]{PT}, which implies
%for each integral with respect constants $\sigma$ and $\overline{\sigma}$ in order to
%which get the following behaviour,
\begin{equation}\label{D^sK:RF:behaviour}
\left| |D|^\theta[K^\beta_\gamma(1,\cdot)](X) \right| \lesssim \frac{1}{|X|^{1+\theta}}, \quad |X| \gg 1.
\end{equation}
Since $\theta>0$, we can apply the lemma if the condition
\[
\sigma, \overline{\sigma} \in \left\{ a+ i b \in \C : \ -\cos\left( \frac{\beta \pi}{2} \right) \leq a \leq 1, \ |b| \leq -\tan\left( \frac{\beta \pi}{2} \right) \right\},
\]
is satisfied. This holds, since  
\[
|\gamma|\leq \min\{\beta,2-\beta\} \ \ \Longrightarrow \ \ \frac{\gamma \pi}{2} \in \left( - \frac{(2-\beta) \pi}{2}, \frac{(2-\beta) \pi}{2}\right) \subset \left( -\frac{\pi}{2}, \frac{\pi}{2}\right)
\]
and this implies that, if $\sigma=a+ib$, 
\[
\cos\left( \frac{(2-\beta) \pi}{2} \right) = - \cos\left( \frac{\beta \pi}{2}\right) \leq  a = \cos\left(\frac{\gamma \pi}{2}\right) \leq 1
\]
and the imaginary part satisfies
\[
|b| = \left| \sin\left(\frac{\gamma \pi}{2}\right) \right| \leq \sin\left( \frac{(2-\beta)\pi}{2} \right) \leq \tan\left( \frac{(2- \beta)\pi}{2} \right) = -\tan\left( \frac{\beta \pi}{2}\right).
\]
As a result of (\ref{D^s:K:RF}) and the previous behaviour given in (\ref{D^sK:RF:behaviour}), we conclude that 
\[
|D|^\theta[K^\beta_\gamma(1,\cdot)](X) \in L^p(\R) \quad \text{for any} \ 1\leq p \leq \infty.
\]
Finally, the fourth estimate follows similarly. we leave the details to the reader.
%In particular, one can conclude a the $L^p$ bound, since \cite[Lemma~2]{PT} can also be applied to obtain
%\[
%\begin{split}
%  | |D|^\theta[\partial_X K(1,\cdot)](X) |
%  % =& \frac{1}{\sqrt{2\pi}} \left| \int_{\R} |\xi|^s (i \xi) \,  e^{-|\xi|^\beta e^{i \sgn(\xi) \gamma \frac{\pi}{2}}} e^{i X \xi} \, d\xi \right| \\ 
%%% =& \frac{1}{\sqrt{2\pi}} \Bigg| i  \int_{0}^{\infty} \xi^{1+s} \, e^{-\xi^{\beta} \left(\cos\left(\frac{\gamma \pi}{2}\right) + i \sin\left(\frac{\gamma \pi}{2}\right)  \right)} e^{i X \xi} \, d\xi \\  
%%% &-  i  \int_{-\infty}^{0}  (-\xi)^{1+s} \, e^{-(-\xi)^{\beta} \left(\cos\left(\frac{\gamma \pi}{2}\right) - i \sin\left(\frac{\gamma \pi}{2}\right)  \right)} e^{i X \xi} \, d\xi \Bigg| \\
%% \leq& \frac{1}{\sqrt{2\pi}} \left| \int_{0}^{\infty} \xi^{1+s} \, e^{-\xi^{\beta} \left(\cos\left(\frac{\gamma \pi}{2}\right) + i \sin\left(\frac{\gamma \pi}{2}\right)  \right)} e^{-i (-X) \xi} \, d\xi \right| \\ 
%% &+  \frac{1}{\sqrt{2\pi}}  \left|\int_{0}^{\infty}  \xi^{1+s} \, e^{-\xi^{\beta} \left(\cos\left(\frac{\gamma \pi}{2}\right) - i \sin\left(\frac{\gamma \pi}{2}\right)  \right)} e^{- i X \xi} \, d\xi \right| \\
%\lesssim \frac{1}{|X|^{2+\theta}}, \quad \text{for} \quad |X| \gg 1.
%\end{split}
%\]

%We conclude the fourth estimate using the previous inequality and the self-similarity property of the kernel as is done in (\ref{D^sK:RF:estimate:t}).
\endproof

Now with Lemma~\ref{K:RF:time:behaviour}, we can proceed as in the proof of \cite[Proposition~3.1]{IS}, to obtain Proposition~\ref{u:p:regularity} for (\ref{RFCL}).

In order to conclude the corresponding weak viscous entropy inequality, similar to Theorem~\ref{weak:entropy:inequality} and the Oleinik type of inequality and all other {\it a priori} estimates, similar to Proposition~\ref{Oleinik:inequality} and Lemma~\ref{properties:u}, for positive solutions, we need the following lemma: 
\begin{lem}[Partial integration by parts and energy estimate]\label{Op:splitting:energy}
Let $\beta \in (1,2)$ and $|\gamma| \leq \min\{ \beta , 2-\beta \}$, then 
\begin{enumerate}
\item[(i)] For functions $g$ and $h$ such that $D^\beta_\gamma[g]$, $\DD^{\theta_1}[g]$, $h$, $\overline{\DD^{\theta_2}}[h]  \in L^2(\R)$, then
\[
 \int_{\R} D^\beta_\gamma[g](x) \, h(x) \, dx =  -\frac{1}{d_{\beta+1}}  (c_\gamma^1+c_\gamma^2) \int_{\R} \DD^{\theta_1} [g](x) \, \overline{\DD^{\theta_2}} [h](x) \, dx
 \]
where $1/2<\theta_1$, $\theta_2<1$ and $\beta = \theta_1 + \theta_2$.
\item[(ii)] 
Moreover, for $1<\beta<2$ and $g$, $D^\beta_\gamma[g] \in L^2(\R)\cap C_b^2(\R)$, we have
\[
- \int_{\R} g(x) \, D^\beta_\gamma[g](x) \, dx \geq 0.
\] 
\end{enumerate}
\end{lem} 
\proof
We note that Lemma~\ref{partial:int:parts:D:alpha} has the easy generalisation %(maybe include in the Lemma...)
\[
%\label{partial:int:parts:adjointD:alpha}
\int_{\R} \partial_x \overline{\DD^{\beta-1}}[g](x) \, h(x) \, dx =  -\int_{\R} \DD^{\theta_1}[g](x) \, \overline{\DD^{\theta_2}}[h](x) \, dx 
\]
for $1/2<\theta_1$, $\theta_2<1$ with $\beta = \theta_1 + \theta_2$. This and Lemma~\ref{partial:int:parts:D:alpha} with $\alpha=\beta-1$ implies (i), using the representation (\ref{RF:integral:representation}).

In order to show (ii), we use again the representation (\ref{RF:integral:representation}) and Lemma~\ref{integrationbyparts}, this gives
\[
\begin{split}
 \int_{\R} g \, D^\beta_\gamma[g](x) \, dx =& \frac{1}{d_{\beta+1}}\left( c^1_\gamma \int_{\R} g \, \partial_x \DD^{\beta-1}[g](x) \, dx + c^2_\gamma \int_{\R} g \, \partial_x \overline{\DD^{\beta-1}}[g](x) \, dx \right)\\
=&\frac{1}{d_{\beta+1}} (c^1_\gamma + c^2_\gamma) \int_{\R} g \, \partial_x \DD^{\beta-1}[g](x) \, dx \leq 0,
\end{split}
\]
where the last inequality is proved as in e.g. \cite{CA}.
\endproof

Part (i) of the above lemma, allows to prove a weak entropy inequality, namely:
\[%begin{equation}
\begin{split}
\int_{0}^{\infty} \int_{\R} \Big( |u(t,x) - k| \partial_t \varphi + \frac{1}{q}\sgn(u(t,x)-k)(|u(t,x)|^{q-1}u(t,x)-|k|^{q-1}k) \partial_x \varphi& \\
+ |u(t,x) - k| \overline{D^\beta_\gamma}[\varphi(t,\cdot)](x)  \Big) \, dx dt  &\geq 0,
\end{split}
%\label{Kruzhkov:entropy:inequality:g}
\]%end{equation}
where, as we have defined also in \cite{DC},
\[
\overline{D^\beta_\gamma}[g](x) = \frac{1}{d_{\beta+1}}\left( c_\gamma^1\partial_x\overline{\DD^{\beta-1}}[g](x)+ c_\gamma^2\partial_x \DD^{\beta-1}[g](x) \right).
\]

We observe that the above lemma is necessary to conclude the analogous of Lemma~\ref{properties:u}, in particular property (vii). Indeed, we need an energy estimate similar to (\ref{energy:estimate:2}). Let us briefly indicate how this is obtained. First, we multiply the equation by $u$ and integrate by parts
 \[
%\begin{equation}\label{energy:estimate:g}
\frac{1}{2}\frac{d}{dt} \int_{\R} u^2 \, dx - \int_{\R} u  D^\beta_\gamma[u](x) \, dx =    0.
%\end{equation}
 \] 
 Now, using (i) above, we obtain the energy type of identity:
 \[
 \frac{1}{2}\frac{d}{dt} \int_{\R} u^2 \, dx+
 \frac{1}{d_{\beta+1}}  (c_\gamma^1+c_\gamma^2) \int_{\R} \DD^{\beta/2} [u](x) \, \overline{\DD^{\beta/2}} [u](x) \, dx=0.
%
% - \int_{\R} u  D^\beta_\gamma[u](x) \, dx =    0.
\]
The second term is positive by (ii), this means that, in fact, 
\[
\frac{1}{2}\frac{d}{dt} \int_{\R} u^2 \, dx+
 \frac{1}{d_{\beta+1}}  (c_\gamma^1+c_\gamma^2) \int_{\R}| \DD^{\beta/2} [u](x) |^2\, dx=0.
 \]
 
The rest of the argument follows unchanged, combining all the results that we have mentioned. Thus, we can generalise the large time asymptotic result Theorem~\ref{subcritical:asymptotics} for the equation (\ref{RFCL}) in the sub-critical case, $1<q<\beta$, for non-negative solutions, obtaining the same rate of convergence.

\bigskip

{\bf Acknowledgements:} This paper was partially supported by the MICINN project PGC-2018-094522-B-I00, and by the Basque Government research group grant IT1247-19. X. Diez-Izagirre also acknowledges the support of the Basque Government through the predoctoral grant PRE-2019-2-0006. The authors would also like to thank Diana Stan for her valuable suggestions.

\bibliographystyle{plain}

\end{document}